\documentclass[12pt,a4paper,twoside]{article}
\usepackage{amssymb}
\begin{document}
\author{S. Albeverio $^{1},$ Sh. A. Ayupov $^{2,}\footnote{ Corresponding author},$ \ \ K. K.
Kudaybergenov  $^3$}
\title{\bf NON  COMMUTATIVE  ARENS ALGEBRAS AND THEIR
DERIVATIONS}

\maketitle

\medskip
$^1$ Institut f\"{u}r Angewandte Mathematik, Universit\"{a}t Bonn,
Wegelerstr. 6, D-53115 Bonn (Germany); SFB 611, BiBoS; CERFIM
(Locarno); Acc. Arch. (USI), \emph{albeverio@uni-bonn.de}

$^2$ Institute of Mathematics, Uzbekistan Academy of Science, F.
Hodjaev str. 29, 700143, Tashkent (Uzbekistan), e-mail:
\emph{sh\_ayupov@mail.ru, e\_ayupov@hotmail.com, mathinst@uzsci.net}

 $^{3}$ Institute of Mathematics, Uzbekistan
Academy of Science, F. Hodjaev str. 29, 700143, Tashkent
(Uzbekistan), e-mail: \emph{karim2006@mail.ru}

\newpage

\begin{abstract}

 Given a von Neumann algebra $M$ with
  a faithful normal semi-finite trace  $\tau,$ we consider the  non
  commutative Arens algebra $L^{\omega}(M, \tau)=\bigcap\limits_{p\geq1}L^{p}(M, \tau)$
  and the related algebras $L^{\omega}_2(M, \tau)=\bigcap\limits_{p\geq2}L^{p}(M, \tau)$ and
  $M+L^{\omega}_2(M, \tau)$ which are proved to be   complete
metrizable locally convex *-algebras. The main purpose of the present paper
is to prove that any  derivation of the algebra $M+L^{\omega}_2(M,
\tau)$ is inner and all derivations of the algebras
$L^{\omega}(M,\tau)$ and $L^{\omega}_2(M, \tau)$ are spatial and
implemented by  elements of $M+L^{\omega}_2(M, \tau).$



\
\end{abstract}
\newpage
\section*{\center 1. Introduction}

The present paper is devoted to  the study of derivations on certain
classes of unbounded operator algebras.

Given a (complex) algebra $\mathcal{A},$ a linear operator  $d:
\mathcal{A}\rightarrow \mathcal{A}$
   is called a  \emph{derivation} if $d(xy)=d(x)y+xd(y)$ for all $x, y\in \mathcal{A}.$
   Each element $a\in \mathcal{A}$ generates a derivation $d_a:\mathcal{A}\rightarrow \mathcal{A}$
   defined as $d_a(x)=ax-xa,\,  x\in \mathcal{A}.$ Such derivations are called  \emph{inner} derivations.

If the element $a$ implementing the derivation $d_a$ belongs to a
larger algebra $\mathcal{B}$ containing $\mathcal{A}$ then $d_a$ is
called a \emph{spatial} derivation.

It is a general algebraic problem to find  algebras which admit
only inner derivations. Such examples are

-- finite dimensional simple central algebras;

-- simple unital $C^{\ast}$-algebras;

-- algebras $B(X)$ of all bounded linear operators on a Banach space
$X$  (cf. [8], [16]).

A more general problem is the following one:  given an algebra
$\mathcal{A},$ does there exist an algebra $\mathcal{B}$ such that

(\emph{i}) $\mathcal{A}$ is an ideal in $\mathcal{B},$ so that any
element $a\in \mathcal{B}$ defines a derivation on $\mathcal{A}$
 by $d_a(x)=ax-xa,\,  x\in \mathcal{A};$

 (\emph{ii}) any derivation of $\mathcal{B}$ is inner;

 (\emph{iii}) any derivation of the algebra $\mathcal{A}$ is spatial and
 implemented by an element from $\mathcal{B}?$

Examples of algebras for which the answer is positive are

-- simple (non unital)  $C^{\ast}$-algebras;

-- the algebra $\mathcal{F}(X)$ of finite rank operators on an
infinite dimensional Banach space $X;$

-- more general standard operator algebras on $X,$ i. e. subalgebras
of $B(X)$ which contain $\mathcal{F}(X)$ (cf. [8], [16]).

 The theory of derivations in operator algebras is an important and well investigated part of the general theory
 of operator algebras, with applications in mathematical physics (see e. g. [7],
 [16], [17]). It is well known that every derivation of a
  $C^{\ast}$-algebra is norm continuous and that every derivation of a von
 Neumann algebra is  inner.
 For  a detailed exposition of the theory of bounded derivations we refer to the monographs of
 Sakai [16], [17]. A comprehensive study of derivations in general
 Banach algebras is given in the monograph of Dales [9] devoted to
 automatic continuity of derivations on various classes of Banach
 algebras.

 Investigations of general unbounded derivations (and derivations on unbounded operator algebras)
 began much later and were motivated mainly by needs of mathematical physics, in particular by the
 problem of constructing the dynamics in quantum statistical mechanics.

 The development of a non commutative integration theory was
  initiated by I. Segal [18],
 who considered new classes of (not necessarily Banach) algebras of unbounded operators, in particular the
 algebra  $L(M)$ of all measurable operators affiliated with a von
 Neumann algebra $M.$ Algebraic, order and topological properties of
 the algebra $L(M)$ are somewhat similar to those of von Neumann algebras, therefore in [4], [5]
 we initiated the study of derivations on the algebra $L(M).$ In the particular commutative case where
 $M=L^{\infty}(0; 1)$ is the algebra of all essentially bounded measurable complex functions on
 $(0; 1),$ the algebra $L(M)$ is isomorphic to the algebra $L^{0}(0; 1)$  of all measurable functions on
 $(0; 1).$ Recent results of  [6] (see also [12]) show that  $L^{0}(0; 1)$  admits non-zero
 (and hence discontinuous) derivations. Therefore the properties of derivations on
  the unbounded operator algebra $L(M)$ are very far from being similar to those on $C^{\ast}$- or von Neumann algebras.

  There are many other classes of unbounded operator algebras, which are important in analysis
  and mathematical physics like     $O^{\ast}$-algebras, $GB^{\ast}$-algebras, $EW^{\ast}$-algebras  (see, e. g., [2], [19]). These algebras  also can be equipped
  by appropriate  topologies and it is natural to study properties
  (such as automatic continuity, innerness, spatiality
  etc) of their derivations. It is known
  that each derivation of the maximal $O^{\ast}$-algebra
  $\mathcal{L}^{+}(\mathcal{D})$ is inner, and if a subalgebra
  $\mathcal{A}$ of $\mathcal{L}^{+}(\mathcal{D})$
  contains all finite rank operators on
   $\mathcal{D},$ then any derivation of $\mathcal{A}$
   is spatial and implemented by an element of
   $\mathcal{L}^{+}(\mathcal{D}),$  where $\mathcal{D}$ is a
  dense linear subspace of a Hilbert space (see for details
  [19, Proposition 6.3.2, Corollary 6.3.3]).

  Interesting examples of the mentioned algebras are given by the Arens algebra
  $L^{\omega}(0,1)=\bigcap\limits_{p\geq1}L^{p}(0,1),$
    introduced in [3], and by its non commutative generalizations
    $L^{\omega}(M, \tau)=\bigcap\limits_{p\geq1}L^{p}(M, \tau),$ where  $M$ is a von Neumann
    algebra with a faithful normal semi-finite trace  $\tau,$
and  $L^{p}(M, \tau)=\{x\in L(M):\tau(|x|^{p})<\infty\}.$
      Non commutative Arens algebras were introduced by Inoue [10] and their properties were investigated in
       [1], [22].

 The main purpose of the present work is to give a complete description of derivations on the non commutative Arens algebras
  $L^{\omega}(M,\tau)$ and related algebras. In particular for
  these algebras we obtain the complete solution of the problems
  mentioned above.

  In Section 2 given a von Neumann algebra $M$ with
  a faithful normal semi-finite trace  $\tau,$ along with the non
  commutative Arens algebra $L^{\omega}(M, \tau)$ we consider some basic
  properties of the related algebras  $L^{\omega}_2(M,
\tau)=\bigcap\limits_{p\geq2}L^{p}(M, \tau)$ and $M+L^{\omega}_2(M,
\tau)$ and prove that they are also complete metrizable locally
convex *-algebras. Applying the theory of Banach pairs from [13] we
describe the predual space of the algebra  $M+L^{\omega}_2(M,
\tau).$ This result enables us to apply  *-weak compactness of
closed bounded sets in the algebra $M+L^{\omega}_2(M, \tau)$ for the
proof of the main results. Namely, in Section 3 we prove that all
derivations on the algebra $M+L^{\omega}_2(M, \tau)$ are inner and
any derivation of the Arens algebra $L^{\omega}(M, \tau)$ is
automatically continuous.

Since the algebras $L^{\omega}(M, \tau)$ and $L^{\omega}_2(M, \tau)$
are (two sided) ideals in $M+L^{\omega}_2(M, \tau)$ any element
$a\in M+L^{\omega}_2(M, \tau)$ defines a derivation on
$L^{\omega}(M, \tau)$ (respectively on $L^{\omega}_2(M, \tau)$) by
$d_a (x)=ax-xa, \, x\in L^{\omega}(M, \tau)$ (respectively $x\in
L^{\omega}_2(M, \tau)$). The main results (Theorem 3.7 and 3.8)
assert that  any derivation of the Arens algebra
$L^{\omega}(M,\tau)$ (respectively of the algebra $L^{\omega}_2(M,
\tau)$) is spatial and implemented by some element $a\in
M+L^{\omega}_2(M, \tau).$ In particular if the trace $\tau$ is
finite then all the above algebras coincide and therefore all
derivations on the Arens algebra  $L^{\omega}(M, \tau)$ are inner.
As a corollary we obtain that commutative Arens algebras (in
particular the algebra $L^{\omega}(0; 1)$) admit only zero
derivations.

\begin{center} {\bf 2. Non commutative Arens algebras}
\end{center}

Let  $M$ be a von Neumann algebra with a faithful normal semi-finite
trace $\tau,$ and  denote by $L(M)$ the algebra of all measurable
operators affiliated with  $M.$

Given $p\geq1$ put $L^{p}(M, \tau)=\{x\in
L(M):\tau(|x|^{p})<\infty\}.$ It is known that  $L^{p}(M, \tau)$ is
a Banach space with respect to the norm
$$\|x\|_p=(\tau(|x|^{p}))^{1/p},\quad x\in L^{p}(M, \tau).$$

   Consider the intersection
\begin{center}
$L^{\omega}(M, \tau)=\bigcap\limits_{p\geq1}L^{p}(M, \tau).$
\end{center}
It is proved in [1] that  $L^{\omega}(M, \tau)$ is a locally convex
complete metrizable  $\ast$-algebra with respect to the  topology
$t$ generated by the family of norms $\{\|\cdot\|_p\}_{p\geq1}.$
Moreover the topology can be defined also by the countable system
(sequence) of norms
$$\|x\|_n'=\max\{\|x\|_1, \|x\|_n\},\,n\in\mathbb{N}.$$

 The algebra
$L^{\omega}(M, \tau)$ is called a (non commutative) \emph{Arens
algebra}.

Non commutative Arens algebras are special cases of
$O^{\ast}$-algebras in the sense of K. Schm\"{u}dgen (see, e. g., [2],
[19]). If the trace $\tau$ is finite then $L^{\omega}(M, \tau)$ is
also an $EW^{\ast}$-algebra [22]. The dual space for $(L^{\omega}(M,
\tau), t)$ was completely described in [1], where it was also proved
that $(L^{\omega}(M, \tau), t)$ is a reflexive space if and only if
the trace $\tau$ is finite.

Recall that a subset $K$ in a linear topological space $E$
 is said to be bounded if given any zero neighborhood $V$ in
$E$ there exists $\alpha>0$ such that $K\subset\beta V$ for all
$\beta>\alpha.$

Since  $(L^{\omega}(M, \tau), t)$ is a countably normed space, a
subset  $K$ in  $L^{\omega}(M, \tau)$ is bounded if and only if
there exists a sequence $(r_n)_{n\in\mathbb{N}}$ of positive numbers
such that $\|x\|_n\leq r_n$ for all $x\in K$ and $n\in\mathbb{N}$
(see  [21], p. 368).

 A linear topological space $E$ is said to be locally bounded if it  admits a bounded neighborhood
 of zero.

Automatic continuity of derivations on Banach algebras and more
general locally bounded  $F$-algebras was investigated in the
monograph [9] of H. G. Dales. Arens algebras are not locally bounded
in general; morever it is not difficult to see that
  $L^{\omega}(M,
\tau)$ is locally bounded if and only if the von Neumann algebra
 $M$ is finite dimensional. Therefore most of results from [9] can not be applied to
 the case of Arens algebras.

Now let us recall the notion of a Banach pair  (see [13]).

Let  $(E, t_E)$ be a Hausdorff linear topological space over the
field of complex numbers, and let $(A, \|\cdot\|_A)$ and $(B,
\|\cdot\|_B)$ be Banach spaces which are linear subspaces of $(E,
t_E)$ such that the topology  $t_E$ induces on $A$ and $B$
topologies which are weaker than the topologies defined by the norms
$\|\cdot\|_A$ and $\|\cdot\|_B,$ respectively. This means exactly
that $A$ and $B$ are topologically imbedded into $(E, t_E).$ In this
case we say that  $A$ and $B$ define a\emph{ Banach pair}. Each
Banach pair defines a couple of Banach spaces $A\cap B$ and $A+B$
with the norms
$$\|x\|_{A\cap B}=\max\{\|x\|_A, \|x\|_B\}, \quad x\in A\cap
B,$$
$$\|x\|_{A+B}=\inf\{\|a\|_A+\|b\|_B: x=a+b,\, a\in A, b\in B\},$$
respectively.

A Banach space  $(Z, \|\cdot\|_Z)$ is said to be \emph{intermediate}
for the Banach pair $A$ and $B,$ if the continuous embeddings
$$A\cap B\subset Z \subset A+B$$ are valid.

Further we shall need the following result from [13, p. 27].

\textbf{Lemma 2.1.} \emph{If the intersection $A\bigcap B$ is dense
in each member of the Banach pair  $A$ and $B,$ then the dual spaces
$A'$ and $B'$ also form a Banach pair. Morever $(A\bigcap B)'$ is
isometrically isomorphic to the space $A'+B',$ and $(A+B)'$ is
isometrically isomorphic to the space $A'\cap B'.$}

 It is known that [20] for the Banach pair  $L^{p_1}(M, \tau)$ and $L^{p_2}(M, \tau), p_1<p_2$
any space  $L^{p}(M, \tau),$ $p\in [p_1; p_2],$ is intermediate.
Therefore
$$\|x\|_p\leq c(p_1,p_2,p) \max \{\|x\|_{p_1}, \|x\|_{p_2}\}$$
for all $x\in L^{p_1}(M, \tau)\cap L^{p_2}(M, \tau),$ where $c(p_1,
p_2,p)$  is a fixed positive number depending on $p_1, p_2, p.$

On the spaces   $$\bigcap\limits_{n=2}^m L^n(M, \tau)
\quad\mbox{and} \quad \sum\limits_{n=2}^m L^{q_n}(M, \tau),$$
consider respectively the norms
$$\|x\|_m^0=\max\{\|x\|_n:
n=\overline{2,m}\},\,\,\,\, x\in \bigcap\limits_{n=2}^m L^n(M,
\tau),$$

$$\|y\|_m'=\inf \{\sum\limits_{n=2}^m\|y_n\|_{q_n}: y_n\in
L^{q_n}(M, \tau), y=\sum\limits_{n=2}^m y_n\},$$
 where
$\frac{\textstyle1}{\textstyle n}+\frac{\textstyle 1}{\textstyle
q_n}=1,$ $n=\overline{2,m},\,m\in\mathbb{N}\setminus\{1\}.$

Consider the Banach pair  $L^{p_1} (M, \tau)$ and $L^{p_2} (M,
\tau),$ $2\leq p_1\leq p_2.$ Since the space
$$M\bigcap L^1(M, \tau)=\{x\in M: \tau(|x|)<\infty\}$$
is dense in all spaces  $L^p(M, \tau),$ $p\geq 1,$ the intersection
$$L^{p_1} (M, \tau)\bigcap L^{p_2} (M, \tau)$$ is  dense in $L^{p_i}
(M, \tau),$ $i=1,2.$  Therefore from Lemma 2.1 by induction on $m$
we obtain the following result:

\textbf{Proposition 2.2. } \emph{The dual space for
$(\bigcap\limits_{n=2}^m L^n(M, \tau) , \|\cdot\|_m^0 )$ is
isometrically isomorphic to $(\sum\limits_{n=2}^m L^{q_n}(M, \tau),
\|\cdot\|_m' )$ and the dual space for $(\sum\limits_{n=2}^m
L^{q_n}(M, \tau),  \|\cdot\|_m' )$ is isometrically isomorphic to
$(\bigcap\limits_{n=2}^m L^n(M, \tau) , \|\cdot\|_m^0).$ The duality
is given by the bilinear form}
$$\langle x,y\rangle=\tau(xy),\quad x\in \bigcap\limits_{n=2}^m
L^n(M, \tau), y\in \sum\limits_{n=2}^m L^{q_n}(M, \tau).$$

Since  $L^{1}(M, \tau)$ and $\sum\limits_{n=2}^m L^{q_n}(M, \tau)$
also form a Banach pair and $L^{1}(M, \tau)\cap M$ is dense in both
of  $L^{1}(M, \tau)$ and $\sum\limits_{n=2}^m L^{q_n}(M, \tau),$
Lemma  2.1 and Proposition  2.2 imply

\textbf{Proposition 2.3. } \emph{The dual space for the Banach space
$L^{1}(M, \tau)\bigcap(\sum\limits_{n=2}^m L^{q_n}(M, \tau))$ is
isometrically isomorphic to  $M+\bigcap\limits_{n=2}^m L^n(M, \tau).$
Moreover, given any  $f\in (L^{1}(M,
\tau)\bigcap(\sum\limits_{n=2}^m L^{q_n}(M, \tau)))'$ there exists a
unique element  $a\in M+\bigcap\limits_{n=2}^m L^n(M, \tau)$ such
that }
$$f(x)=\tau(xa),\quad x\in L^{1}(M,
\tau)\bigcap(\sum\limits_{n=2}^m L^{q_n}(M, \tau)).$$

Now consider the following space
$$L^{\omega}_2(M, \tau)=\bigcap\limits_{p\geq 2}L^p(M, \tau)$$
with the topology $t_2$ generated by the family of norms
$\{\|\cdot\|_p\}_{p\geq2}.$

\textbf{Proposition 2.4.} \emph{$(L^{\omega}_2(M, \tau), t_2)$ is a
complete metrizable locally convex  *-algebra.}

 Proof. From the inequality $\|xy\|_p\leq\|x\|_{2p}\|y\|_{2p}$ it easily
 follows that $L^{\omega}_2(M, \tau)$ is closed under
 the multiplication. It is also clear that
$L^{\omega}_2(M, \tau)$ is closed under the involution, i. e. it
forms a $\ast$-algebra.

Let us show that

$$\bigcap\limits_{n=2}^{\infty} L^{n}(M, \tau)=L_2^{\omega}(M, \tau).$$

Clearly  $L_2^{\omega}(M, \tau)\subset
 \bigcap\limits_{n=2}^{\infty} L^{n}(M, \tau).$
For any $p\geq2$ take natural numbers $n_1$ and $n_2$ such that
$n_1\leq p\leq n_2.$ Since the space  $L^p(M, \tau)$ is intermediate
for the Banach pair $L^{n_1}(M_1, \tau)$ and $L^{n_2}(M_2, \tau),$
we have
$$ L^{n_1}(M, \tau)\cap L^{n_2}(M, \tau)\subset
L^{p}(M, \tau)$$ and
$$\|x\|_p\leq c(n_1 , n_2 , p)\max\{\|x\|_{n_1}, \|x\|_{n_2}\}\leq
c(n_1 , n_2 , p)\|x\|_{n_2}^{0}.$$ This means that
$\bigcap\limits_{n=2}^{\infty} L^{n}(M, \tau)=L_2^{\omega}(M,
\tau),$  and the topology $t_2$ is generated by the system of norms
 $\{\|\cdot\|_{n}^{0}\}_{n\geq2} .$

Let us show that  $(L^{\omega}_2(M, \tau), t_2)$ is complete. Let
$(x_n)_{n\in \mathbb{N}}$  be a Cauchy sequence in $(L^{\omega}_2(M,
\tau), t_2).$ Then $(x_n)_{n\in \mathbb{N}}$ is a Cauchy sequence
 in $(L^p(M, \tau),\\\|\cdot\|_p)$ for all $p\geq 2$ and
hence there exists  $a_p\in L^p(M, \tau)$ such that
$\|x_n-a_p\|_p\rightarrow 0$  as $n\rightarrow\infty.$ We have
$a_p=a_q$ for all $p, q\geq 2.$ Indeed, let $e\in M$ be a projection
and $\tau(e)<\infty.$ Then
$$(a_p-a_q)e\in L^2(M, \tau),$$
moreover
$$\|(a_p-a_q)e\|_2\leq \|(a_p-x_n)e\|_2+\|(x_n-a_q)e\|_2\leq
\|e\|_{p_1}\|(a_p-x_n)\|_p+\|e\|_{q_1}\|\|(x_n-a_q)\|_q,$$ where
$p_1, q_1\in(2; \infty],$ $\frac{\textstyle
 1}{\textstyle  p_1}+\frac{\textstyle  1}{\textstyle p}=\frac{\textstyle 1}{\textstyle 2},$ $\frac{\textstyle 1}{\textstyle
q_1}+\frac{\textstyle
 1}{\textstyle  q}=\frac{\textstyle  1}{\textstyle  2}.$
Therefore  $\|(a_p-a_q)e\|_2=0$ and $a_p e=a_q e$ for each
projection $e\in M$ with $\tau(e)<\infty.$ Since $\tau$ is a
semi-finite trace this means that $a_p=a_q.$

Therefore, $L^{\omega}_2(M, \tau)$ is a locally convex complete
metrizable  *-algebra, with respect to  the topology $t_2,$
generated by the family of  norms $\{\|\cdot\|_n\}_{n\geq 2}.$ The proof is complete. $\blacksquare$

Note that if $\tau(\textbf{1})<\infty$ then $L^{\omega}_2(M,
\tau)=L^{\omega}(M, \tau),$ and the topology  $t_2$ coincides with
the topology $t.$

 On the space  $M+
L^{\omega}_2(M, \tau)$ consider the family of  norms
$\{\|\cdot\|''_{n}\}_{n\geq2}$ defined by
$$\|x\|''_{n}=\inf\{\|x_1\|_{\infty}+\|x_2\|_{n}^{0}: x=x_1+x_2,\, x_1\in M,\, x_2\in
L^{\omega}_2(M, \tau)\}.$$

Let  $t_0$ be  the topology on  $M+ L^{\omega}_2(M, \tau)$ generated
by the family of norms  $\{\|\cdot\|''_{n}\}_{n\geq2}.$

\textbf{Lemma 2.5.} \emph{Let  $(x_n)_{n\in\mathbb{N}}$ be a
sequence in $M+ L^{\omega}_2(M, \tau)$ such  that
$\|x_n\|_{n}''\rightarrow0$ as $n\rightarrow\infty.$ Then}
$x_n\stackrel{t_0}{\longrightarrow}0$ \emph{as}
$n\rightarrow\infty.$

Proof. Let  $k\in \mathbb{N}\setminus\{1\}.$ Since
$\|x_n\|_{k}''\leq\|x_n\|_{n}''$ for $n\geq k,$ then
 $\|x_n\|_{k}''\rightarrow0$ as $n\rightarrow\infty$ for all  $k\geq 2.$
 Therefore $x_n\stackrel{t_0}{\longrightarrow}0$ as
 $n\rightarrow\infty.$ The proof is complete. $\blacksquare$

\textbf{Proposition 2.6.}  \emph{The algebra ($M+ L^{\omega}_2(M,
\tau), t_0)$ is a complete metrizable locally convex $\ast$-algebra.
Moreover $L^{\omega}(M, \tau)$ is an ideal in $M+L^{\omega}_2(M,
\tau).$ }

Proof. For $x\in M$ and $y\in L^{\omega}_2(M, \tau)$ it is clear
that $xy, yx\in L^{\omega}_2(M, \tau).$ Since $L^{\omega}_2(M,
\tau)$ is closed  under the multiplication, it follows that
$M+L^{\omega}_2(M, \tau)$  forms a $\ast$-algebra.

Take $x, y \in M+ L^{\omega}_2(M, \tau).$ We have
$$ \|xy\|_{n}''\leq\|x\|_{2n}''\|y\|_{2n}'' .$$

Indeed, let  $\varepsilon>0.$ Take   $x_1 , y_1 \in M,$ $x_2, y_2
\in L^{\omega}_2(M, \tau)$ such that $x=x_1+x_2 ,$ $y=y_1+y_2 ,$ $
\|x\|_{n}''\geq\|x_1\|_{\infty}+\|x_2\|_{n}^{0}-\varepsilon,$ $
\|x\|_{n}''\geq\|y_1\|_{\infty}+\|y_2\|_{n}^{0}-\varepsilon.$

Then $ \|xy\|_{n}''= \|x_1 y_1+x_1 y_2+x_2 y_1+x_2
y_2\|_{n}''\leq\|x_1 y_1\|_{\infty}+\|x_1 y_2\|_{n}^{0}+
 \|x_2 y_1\|_{n}^{0}+\|x_2 y_2\|_{n}^{0}=\|x_1 y_1\|_{\infty}+
 \max\limits_{2\leq i\leq n}\{\|x_1 y_2\|_{i}\}+
 \max\limits_{2\leq i\leq n}\{\|x_2 y_1\|_{i}\}+
 \max\limits_{2\leq i\leq n}\{\|x_2 y_2\|_{i}\}\leq\|x_1 y_1\|_{\infty}+
 \|x_1\|_{\infty}\max\limits_{2\leq i\leq n}\{\|y_2\|_{i}\}+
\|y_1\|_{\infty} \max\limits_{2\leq i\leq n}\{\|x_2\|_{i}\}+
 \max\limits_{2\leq i\leq n}\{\|x_2\|_{2i}\}\max\limits_{2\leq i\leq n}\{\|y_2\|_{2i}\}
 \leq\|x_1\|_{\infty}\| y_1\|_{\infty}+
 \|x_1\|_{\infty}\|y_2\|_{2n}^{0}+
\|y_1\|_{\infty}\|x_2\|_{2n}^{0}+
\|x_2\|\|_{2n}^{0}\|y_2\|_{2n}^{0}\leq(\|x_1\|_{\infty}
+\|x_2\|_{2n}^{0})(\|y_1\|_{\infty}
+\|y_2\|_{2n}^{0})\leq\\\leq(\|x\|_{2n}''+\varepsilon)(\|y\|_{2n}''+\varepsilon),$
i. e. $
\|xy\|_{n}''\leq(\|x\|_{2n}''+\varepsilon)(\|y\|_{2n}''+\varepsilon).$
Since  $\varepsilon>0$ is arbitrary this implies that  $
\|xy\|_{n}''\leq\|x\|_{2n}''\|y\|_{2n}'' .$

Clearly $ \|x^{\ast}\|_{n}''=\|x\|_{n}''.$

This means that the multiplication and the involution on $M+
L^{\omega}_2(M, \tau)$ are continuous.

Let us show that $M+ L^{\omega}_2(M, \tau)$ is  complete. Let
$(x_n)_{n\in \mathbb{N}}$ be a Cauchy sequence in $M+
L^{\omega}_2(M, \tau).$ Then $\|x_n - x_m\|_{k}''\rightarrow0$ as
$n, m\rightarrow\infty$ for all $k\geq2.$ Therefore, there exists a
subsequence $(x_{n_{k}})_{k\in\mathbb{N}}$ such that
$$\|x_{n_{k+1}} - x_{n_{k}}\|_{k}''<2^{-k-1}, k\geq2.$$
Set  $y_k =x_{n_{k+1}} - x_{n_{k}}, k\geq2.$ Take  $b_k\in M, b_k
'\in L^{\omega}_2(M, \tau)$ such that  $b_k +b_k '=y_k$ and
$$\|b_{k}\|_{\infty}+\|b_k '\|_{k}^{0}<\|y_{k}\|_{k}''+2^{-k-1}<2^{-k}, k\geq2.$$
Then the series $\sum\limits_{k=2}^{\infty}b_k$ and
$\sum\limits_{k=2}^{\infty}b_k '$  converge in $M$ and
$L^{\omega}_2(M, \tau)$ respectively. Put
$b=\sum\limits_{k=2}^{\infty}b_k$ and
$b'=\sum\limits_{k=2}^{\infty}b_k '$ and consider the sums
$s_n=\sum\limits_{k=2}^{n}b_k$ and $s'_n=\sum\limits_{k=2}^{n}b_k
'.$ Then  $\|s_n+s'_n-(b+b')\|_{n}''\rightarrow0$ as
$n\rightarrow\infty,$ and by Lemma  2.5 we have
$s_n+s'_n-(b+b')\stackrel{t_0}{\longrightarrow}0.$ On the other hand
$s_k+s'_k=x_{n_{k+1}}-x_{n_{2}},$ i. e.
$x_{n_{k}}\stackrel{t_0}{\longrightarrow} x_{n_{2}}+b+b'.$ Since
$(x_n)$ is a Cauchy   sequence in  $M+L^{\omega}_2(M, \tau)$ this
implies that $x_n\stackrel{t_0}{\longrightarrow} x_{n_{2}}+b+b'.$

Now we shall  show that  $L^{\omega}(M, \tau)$ is an ideal in
$M+L^{\omega}_2(M, \tau).$

Take  $x\in L^{\omega}(M, \tau)$ and $a\in L^{\omega}_2(M, \tau).$
Then $x\in L^{2p}(M, \tau)$ and $a\in L^{2p}(M, \tau)$ for all
$p\geq1.$ Thus $ax, xa\in L^{p}(M, \tau)$  for all $p\geq1,$ and
therefore  $ax, xa\in L^{\omega}(M, \tau),$ i. e. the algebra
$L^{\omega}(M, \tau)$ is an ideal in $L^{\omega}_2(M, \tau).$ Since
for $x\in L^{\omega}(M, \tau)$ and $a\in M$ we have that $ax, xa\in
L^{\omega}(M, \tau),$ this implies that $L^{\omega}(M, \tau)$ is an
ideal in $M+L^{\omega}_2(M, \tau).$ The proof is complete. $\blacksquare$

\textbf{Remark 1.} In a similar way it follows that the algebra
$L^{\omega}_2(M,\tau)$ is also an ideal in $M+L^{\omega}_2(M,
\tau).$

\textbf{Remark 2.} Note that if $\tau(\textbf{1})<\infty$ then
$M+L^{\omega}_2(M, \tau)=L^{\omega}(M, \tau),$ and the  topology
$t_0$ coincides with the topology $t.$

\textbf{Lemma  2.7.}  $M+ L^{\omega}_2(M,
\tau)=\bigcap\limits_{m=2}^{\infty}(M+
\bigcap\limits_{n=2}^{m}L^{n}(M, \tau)).$

Proof. It is sufficient to  show that
$\bigcap\limits_{m=2}^{\infty}(M+ \bigcap\limits_{n=2}^{m}L^{n}(M,
\tau))\subseteq M+ L^{\omega}_2(M, \tau).$ Take  $a\in
\bigcap\limits_{m=2}^{\infty}(M+ \bigcap\limits_{n=2}^{m}L^{n}(M,
\tau)).$ Then $a\in M+ \bigcap\limits_{n=2}^{m}L^{n}(M, \tau)$ for
all $m\geq 2.$ Therefore there exist  $b_m\in M, c_m\in
\bigcap\limits_{n=2}^{m}L^{n}(M, \tau)$ such that $a=b_m+c_m.$ Since
$L^{1}(M, \tau)\bigcap M$ is dense in
$\bigcap\limits_{n=2}^{m}L^{n}(M, \tau),$ there exists $d_m\in
L^{1}(M, \tau)\bigcap M$ such that $\|c_m-d_m\|_{m}^{0}<1/m$ for all
$ m\geq2.$ Then
$\|a-b_m-d_m\|_{m}''=\|c_m-d_m\|_{m}''\leq\|c_m-d_m\|_{m}^{0}\rightarrow0.$
By Lemma 2.5  we have $b_m+d_m\stackrel{t_0}{\longrightarrow} a.$
Since $M+ L^{\omega}_2(M, \tau)$ is  $t_2$-complete and $b_m+d_m\in
M+ L^{\omega}_2(M, \tau),$ one has $a\in M+ L^{\omega}_2(M, \tau).$ The proof is complete.
$\blacksquare$

Now let us prove the following equality
$$\bigcup\limits_{m=2}^{\infty}\sum\limits_{n=2}^{m} L^{q_n}(M,
\tau)=\mbox{Lin}(\bigcup\limits_{1<q\leq 2}L^{q}(M, \tau)),$$ where
$\frac{\textstyle 1}{\textstyle n}+\frac{\textstyle
 1}{\textstyle  q_n}=1$ and $\mbox{Lin}(E)$ denotes the linear span
 of $E.$

 Let $q\in (1,
2].$  Take $q_n$ and $q_m$ such that $q_n\leq q \leq q_m .$ Since
$L^{q}(M, \tau)$ is intermediate for the Banach pair $L^{q_n}(M,
\tau),$ $L^{q_n}(M, \tau),$ we have
$$L^{q}(M,
\tau)\subset L^{q_n}(M, \tau)+L^{q_n}(M, \tau).$$ Therefore,
$\bigcup\limits_{m=2}^{\infty}\sum\limits_{n=2}^{m} L^{q_n}(M,
\tau)=\mbox{Lin}(\bigcup\limits_{1<q\leq 2}L^{q}(M, \tau)).$

On the space $\bigcup\limits_{m=2}^{\infty}\sum\limits_{n=2}^{m}
L^{q_n}(M, \tau)$  consider the norm

 $$\|y\|_1'=\inf
\{\sum\limits_{n=2}^m\|y_n\|_{q_n}: y_n\in L^{q_n}(M, \tau),
y=\sum\limits_{n=2}^m y_n\},$$
 where
$\frac{\textstyle 1}{\textstyle n}+\frac{\textstyle 1}{\textstyle
q_n}=1,$ $n=\overline{2,m}, m\in \mathbb{N}\setminus\{1\},$\\ and on
the space $L^{1}(M, \tau)\bigcap
(\bigcup\limits_{m=2}^{\infty}\sum\limits_{n=2}^{m} L^{q_n}(M,
\tau))$  define the norm  as

$$\|x\|_1^{0}=\max
\{\|x\|_{1}, \|x\|_{1}'\}, \quad x\in L^{1}(M, \tau)\bigcap
(\bigcup\limits_{m=2}^{\infty}\sum\limits_{n=2}^{m} L^{q_n}(M,
\tau)).$$

The main result of this section is the following theorem which
describes the predual space of the algebra $M+ L^{\omega}_2(M,
\tau).$

 \textbf{Theorem  2.8.} \emph{The dual space for} $L^{1}(M,
\tau)\bigcap(\mbox{Lin}(\bigcup\limits_{1<q\leq2}L^{q}(M, \tau)))$
\emph{is isomorphic to} $M+ L^{\omega}_2(M, \tau).$

Proof. Let  $a=b+c\in M+ L^{\omega}_2(M, \tau).$ Then putting

$$f_a (x)=\tau(xb)+\tau(xc)\eqno (1)$$
we define a continuous linear functional $f_a $ on $L^{1}(M,
\tau)\bigcap(\mbox{Lin}(\bigcup\limits_{1<q\leq2}L^{q}(M, \tau))).$

Conversely, let us show that any continuous linear functional on
\\ $L^{1}(M, \tau)\bigcap(\mbox{Lin}(\bigcup\limits_{1<q\leq2}L^{q}(M, \tau)))$
has the form  (1).

Let $f\in (L^{1}(M,
\tau)\bigcap(\mbox{Lin}(\bigcup\limits_{1<q\leq2}L^{q}(M,
\tau))))'.$ Since the restriction of the   linear functional $f$ on
$L^{1}(M, \tau)\bigcap(\sum\limits_{n=2}^{m} L^{q_n}(M, \tau))$ is
continuous,  by Proposition 2.3 there exists $a_m =b_m+c_m\in
M+\bigcap\limits_{n=2}^{m} L^{n}(M, \tau)$ such that
$$f (x)=\tau(xb_m)+\tau(xc_m),\quad x\in L^{1}(M, \tau)\bigcap M.$$
Thus $\tau(xb_i)+\tau(xc_i)=\tau(xb_j)+\tau(xc_j)$ for all
$i,j\geq2.$ Therefore,
$$\tau(x(b_i-b_j))=\tau(x(c_j-c_i)),\quad x\in L^{1}(M, \tau)\bigcap M.\eqno (2)$$
Since $L^{1}(M, \tau)\bigcap M$ is dense in $L^{1}(M,
\tau)\bigcap(\sum\limits_{n=2}^{m} L^{q_n}(M, \tau))$ for all
$m\in\mathbb{N}\setminus\{1\},$  by (2) we obtain that
$b_i-b_j=c_j-c_i,$ i. e. $a_i=b_i+c_i=b_j+c_j=a_j$ for all
$i,j\geq2.$ Thus $a_i=a_j$ for all $i,j\geq2.$ By Lemma  2.7 it
follows that $a_2\in M+ L^{\omega}_2(M, \tau).$ Since $L^{1}(M,
\tau)\bigcap M$ is dense in $L^{1}(M,
\tau)\bigcap(\mbox{Lin}(\bigcup\limits_{1<q\leq2}L^{q}(M, \tau))),$
we have $f(x)=\tau(xb_1)+\tau(xc_1)$ for all  $x\in L^{1}(M,
\tau)\bigcap(\mbox{Lin}(\bigcup\limits_{1<q\leq2}L^{q}(M, \tau))).$
The uniqueness of the element $a_2$  follows from the density of the
set $L^{1}(M, \tau)\bigcap M$ in $L^{1}(M,
\tau)\bigcap(\mbox{Lin}(\bigcup\limits_{1<q\leq2}L^{q}(M, \tau))).$

Now in a standard way one  proves that the map $a\mapsto f_a$ is a
linear isomorphism between $M+ L^{\omega}_2(M, \tau)$ and $(L^{1}(M,
\tau)\bigcap(\mbox{Lin}(\bigcup\limits_{1<q\leq2}L^{q}(M,
\tau))))'.$ The proof is complete. $\blacksquare$

Proposition 2.8 implies the following

 \textbf{Corollary  2.9.} \emph{If
$\tau(\textbf{1})<\infty$ then the dual space for}
$\mbox{Lin}(\bigcup\limits_{1<q\leq2}L^{q}(M, \tau))$ \emph{is
isomorphic to} $L^{\omega}(M, \tau).$

For general Arens algebras the predual space does not exist but the
dual space was described in  [1] in the following way.

\textbf{Proposition  2.10.} \emph{The dual space for}
$L^{\omega}(M, \tau)$ \emph{is
isomorphic to} $\mbox{Lin}(\bigcup\limits_{1<q\leq\infty}L^{q}(M, \tau)).$

 Since any closed bounded subset in the dual space of a locally convex
 space is $\ast$-weakly compact, Theorem 2.8  implies

\textbf{Corollary 2. 11.} \emph{Each closed bounded subset in $M+
L^{\omega}_2(M, \tau)$ is  $\ast$-weakly compact. }

\begin{center}
{\bf 3. Derivations of Arens algebras}
\end{center}

Let  $\mathcal{A}$ be a complex algebra and  let $E$ be a complex
linear space. Recall that $E$ is called a left $\mathcal{A}$-module
(respectively right $\mathcal{A}$-module) if a bilinear map $(a,
x)\mapsto a\cdot x$ (respectively $(a, x)\mapsto x\cdot a$) from
$\mathcal{A}\times E$ into $E$ is defined
 such that   $$a\cdot(b\cdot
x)=ab\cdot x \quad (\mbox{respectively}\, (x\cdot a)\cdot b=x\cdot
ab),$$ for all $a, b\in \mathcal{A},\, x\in E$.

$E$ is said to be   $\mathcal{A}$-bimodule if  $E$ is simultaneously
a left and right
 $\mathcal{A}$-module such that
 $$a\cdot(x\cdot b)=(a\cdot x)\cdot b,$$ for all  $a, b\in
\mathcal{A}, x\in E.$

 Let  $\mathcal{A}$ be a Banach algebra and let $E$ be a Banach space.

If $E$ is a $\mathcal{A}$-bimodule and the maps  $(a, x)\mapsto
a\cdot x$ and $(a, x)\mapsto x\cdot a$ are continuous, then $E$ is
called a Banach  $\mathcal{A}$-bimodule.

For example, the Banach space  $L^{p}(M, \tau),$ $p\geq1$ is a
Banach  $M$-bimodule. Indeed, since for all $a\in M$ and $x\in
L^{p}(M, \tau)$ one has   $ax\in L^{p}(M, \tau),$ $xa\in L^{p}(M,
\tau)$ and $ \|a x\|_{p}\leq \|a\|_{\infty}\|x\|_p,$ where
$\|\cdot\|_{\infty}$ is the $C^{\ast}$-norm on  $M,$ the space
$L^{p}(M, \tau)$ is a Banach  $M$-bimodule. Therefore, the space
$M+L^{p}(M, \tau)$ ($p\geq1$)  is  also a Banach $M$-bimodule.

Further we need the following result due Ringrose  (see  [15],
Theorem 2, also [9], p. 638).

\textbf{Theorem  3.1}. \emph{Let  $\mathcal{A}$ be a
$C^{\ast}$-algebra and $E$ be a Banach $\mathcal{A}$-bimodule. Then
each derivation
  $D:\mathcal{A}\rightarrow E$  is continuous.}

One of the main results of the present work is the following

 {\bf Theorem 3.2.} \emph{Let  $M$ be
a von Neumann algebra with a faithful normal semi-finite trace
$\tau.$ Then each derivation of the algebra
  $M+L^{\omega}_{2}(M, \tau)$ is inner.}

Proof.  Since $M+\bigcap\limits_{n=2}^{m}L^{n}(M, \tau)$ ($m\geq2$)
is a Banach  $M$-bimodule,   by Theorem 3.1 the derivation $d$  is a
continuous map from  $M$ into $M+\bigcap\limits_{n=2}^{m}L^{n}(M,
\tau)$ for all $m\geq2.$ Therefore,   $\|x_k\|_{\infty}\rightarrow
0$ ($k\rightarrow\infty$) implies $\|d(x_k)\|_{n}''\rightarrow 0$
($k\rightarrow\infty$) for all $n\geq 2.$ Since $M+L^{\omega}_{2}(M,
\tau)$ is a countably normed space with the family of norms
$\{\|\cdot\|_{n}''\}_{n\in \mathbb{N}\setminus\{1\}},$ the operator
$d:M\rightarrow M+L^{\omega}_{2}(M, \tau)$ is  continuous.

Let $U$ be the group of all unitary elements in $M.$ For  $u\in U$
put
$$T_u(x)=uxu^{\ast}+d(u)u^{\ast}, \quad x\in M+L^{\omega}_{2}(M, \tau).$$
Then for $u, v\in U,$ one has
\begin{center}
$T_u(T_v(x))=T_u(vxv^{\ast}+d(v)v^{\ast})=u(vxv^{\ast}+d(v)v^{\ast})u^{\ast}+d(u)u^{\ast}=
uvxv^{\ast}u^{\ast}+ud(v)v^{\ast}u^{\ast}+d(u)u^{\ast}=(uvx+d(u)v+ud(v))(uv)^{\ast}=
uvx(uv)^{\ast}+d(uv)(uv)^{\ast}=T_{uv}(x),$
\end{center}
i. e. $$ T_uT_v=T_{uv}, \quad u, v\in U.\eqno (3)$$

Since the operator  $d$ is  continuous there exist
$C_n>0,\,n\in\mathbb{N}\setminus\{1\},$ such that
$$
\|d(x)\|_n''\leq C_n\|x\|_{\infty},\quad  x\in M.
$$
From  $\|u\|_{\infty}=1 \, (u\in U)$ it follows
$$\|T_u(0)\|_n
''=\|d(u)u^{\ast}\|_n''\leq
\|d(u)\|_{n}''\|u^{\ast}\|_{\infty}=\|d(u)\|_{n}''\leq C_n,$$ i. e.
$\|T_u(0)\|_n ''\leq C_n$ for all  $u\in U.$ Therefore, the set
$K_d=\{T_v(0): v\in U\}$ is bounded. Moreover, the set
$K=\mbox{cl(co}(K_d))$ -- the  closure of the convex hull of $K_d,$
is a closed convex bounded subset in $M+L^{\omega}_2(M, \tau).$ By
Corollary 2.11 $K$ is a non-empty convex *-weakly compact set. By
(3) we have  $T_u(K_d)\subset K_d$ for all $u\in U.$ Since $T_u$ is
an affine homeomorphism,
$T_u(\mbox{cl(co}(K_d)))=\mbox{cl(co}(T_u(K_d)))\subset
\mbox{cl(co}(K_d)),$ ($u\in U$), i. e. $T_u (K)\subseteq K.$

For  $x, y\in M+ L^{\omega}_2(M, \tau)$ we have
$$\|(T_u(x)-T_u(y))\|_{n}''=\|u(x-y)u^{\ast}\|_n
''=\|x-y\|_n ''.$$ Therefore by Ryll-Nardzewski's fixed point
theorem  [14] there exists  $a\in K$ such that  $T_u (a)=a$ for all
$u\in U.$ Therefore $uau^{\ast}+d(u)u^{\ast}=a,$ i. e. $d(u)=au-ua$
for all $u\in U.$ Since any element of $M$ is a finite linear
combination of unitary  elements in $M$, we have $d(x)=ax-xa$ for
all $x\in M.$

Now let us show that  $d(x)=ax-xa$ for all $x\in M+ L^{\omega}_2(M,
\tau).$

First suppose that  $x\in M+ L^{\omega}_2(M, \tau),\, x\geq0.$ Then
the element $\textbf{1}+x$ is  invertible and
$(\textbf{1}+x)^{-1}\in M.$

Let   $b\geq0$ be an invertible elements of $M+L^{\omega}_2(M,
\tau).$ Since  $\textbf{1}=\textbf{1}^{2},$ then
$d(\textbf{1})=d(\textbf{1}^{2})=d(\textbf{1})\textbf{1}+\textbf{1}d(\textbf{1})=2d(\textbf{1}),$
i. e. $d(\textbf{1})=0.$ Therefore
$0=d(\textbf{1})=d(bb^{-1})=d(b)b^{-1}+bd(b^{-1}),$ i. e.
$d(b)=-bd(b^{-1})b.$

Using this equality we obtain
$$d(x)=d(\textbf{1}+x)=-(\textbf{1}+x)d((\textbf{1}+x)^{-1})(\textbf{1}+x).$$
 On  the other hand, since  $(\textbf{1}+x)^{-1}\in M$ one has
$$d((\textbf{1}+x)^{-1})=a(\textbf{1}+x)^{-1}-(\textbf{1}+x)^{-1}a.$$
Therefore, $-(\textbf{1}+x)d((\textbf{1}+x)^{-1})(\textbf{1}+x)=
-(\textbf{1}+x)[a(\textbf{1}+x)^{-1}-(\textbf{1}+x)^{-1}a](\textbf{1}+x)=-(\textbf{1}+x)a+a(\textbf{1}+x)=
ax-xa,$ i. e.
$$d(x)=-(\textbf{1}+x)d((\textbf{1}+x)^{-1})(\textbf{1}+x)=ax-xa.$$
Since any element of $M+L^{\omega}_2(M, \tau)$ is a finite linear
combination of positive  elements in $M+L^{\omega}_2(M, \tau),$ we
have $d(x)=ax-xa$ for all $M+L^{\omega}_2(M, \tau).$  The proof is complete. $\blacksquare$

 {\bf
Corollary 3.3.} \emph{Consider   the Arens
algebra  $L^{\omega}(M, \tau),$ where  $\tau$ is a finite trace. Then any derivation $d$ on
$L^{\omega}(M, \tau)$ is inner. In particular,  $d$ is
$t$-continuous and *-weakly continuous. Moreover, the element $a\in
L^{\omega}(M, \tau)$ implementing $d$ can be taken such that
$\|a\|''_n\leq \|d\|_n, n\geq2.$ }

The continuity of derivations can be proved also in the general case.
Namely, we have the following

\textbf{Proposition  3.4.} \emph{Let  $M$ be a von Neumann algebra
with a faithful normal semi-finite trace $\tau.$ Then each
derivation of the algebra
  $L^{\omega}(M, \tau)$ (respectively  $L^{\omega}_2(M, \tau)$  ) is
   $t$-continuous (respectively $t_2$-continuous).}

Proof. Let us prove the assertion for the algebra $L^{\omega}(M,
\tau),$ the case of $L^{\omega}_2(M, \tau)$ is similar.

 Let $d$
be a derivation on $L^{\omega}(M, \tau).$ Consider the sequence
$(a_n)_{n\in\mathbb{N}}$ in $L^{\omega}(M, \tau),$ such that
$a_n\stackrel{t}\rightarrow 0$ and $d(a_n)\stackrel{t}\rightarrow y$
for some $y\in L^{\omega}(M, \tau).$ Since $L^{\omega}(M, \tau)$ is
a complete metrizable space, by the closed graph theorem it is
sufficient to  show that $y=0.$

Let  $e\in M$ be a projection with finite trace and consider the
Arens algebra $L^{\omega}(eMe, \tau_e)$ associated with the von
Neumannn algebra  $eMe$ and the  faithful normal finite trace
$\tau_e,$ where $\tau_e$ is the restriction of $\tau$ on $eMe.$

Put
$$d_e(x)=ed(exe)e, \quad x\in L^{\omega}(eMe, \tau_e).\eqno (4)$$

For   $x, y\in eMe,$ since $x=exe, y=eye,$ one has
$d_e(xy)=ed(exye)e=ed(exeeye)e=ed(exe)eye+exed(eye)e=d_e(x)y+xd_e(y),$
i. e. $d_e$ is a derivation on $L^{\omega}(eMe, \tau_e).$

 By Corollary 3.3 the derivation  $d_e$ is continuous.
Thus from $ea_n e\stackrel{t}\rightarrow 0$
it follows that $d_e(ea_n e)\stackrel{t}\rightarrow 0$ as
$n\rightarrow\infty.$ On other hand $d_e(ea_n e)=ed(ea_n
e)e=ed(e)a_n e+ed(a_n)e+ea_nd(e)e\stackrel{t}\rightarrow eye$ as
$n\rightarrow\infty.$ Therefore, $$eye=0 \eqno (5)$$ for all
projections  $e\in M$ with finite trace. Since $a_n
y^{\ast}\stackrel{t}\rightarrow0$ and $d(a_n
y^{\ast})=d(a_n)y^{\ast}+a_nd(y^{\ast})\stackrel{t}\rightarrow
yy^{\ast},$ so in (5) the    element $y$ can be a replaced by
$yy^{\ast}.$ Thus $eyy^{\ast}e=0,$ i. e. $(ey)(ey^{\ast})=0,$ and
therefore, $ey=0,$ i. e.
$$y^{\ast}ey=0 \eqno (6)$$ for each
projection  $e\in M$ with finite trace. Since the map $x\mapsto
y^{\ast}xy$ is positive and monotone continuous, taking
$e\uparrow\textbf{1}$ in (6), we obtain  that $yy^{\ast}=0.$
Therefore $y=0.$ The proof is complete. $\blacksquare$

It is easy to see that in commutative Arens algebras any derivation is
equal to zero on projections. Since the linear span of projection is
$t$-dense in any Arens algebra Proposition 3.4 implies

 {\bf Corollary 3.5.} \emph{If  $M$ is  an abelian   von
Neumann algebra  with a faithful normal semi-finite trace $\tau$
then all derivations on
  $L^{\omega}(M, \tau)$ are identically zero.}

\textbf{Remark 3.} As it was noted above the commutative algebra
$L^{0}(0,1)$ of all complex measurable functions on  $(0,1)$ admits
nonzero derivations  (see [6], [12]). On the other hand the
Corollary 3.5 shows that the Arens algebra $L^{\omega}(0,1)$ admits
only zero derivations (similar to the algebra $L^{\infty}(0,1)$),
though it contains  unbounded elements.

The following proposition gives one more type of continuity for
derivations of Arens algebras.

\textbf{Proposition 3.6.} \emph{Let  $d:L^{\omega}(M,
\tau)\rightarrow L^{\omega}(M, \tau)$ be a derivation. Then
 $d$ maps any weakly converging net from $L^{\omega}(M,
\tau)$   into a  net  converging in the *-weak topology in
$M+L^{\omega}_2(M, \tau).$}

Proof. Let a net $(x_\alpha)_{\alpha\in A}\subset L^{\omega}(M,
\tau)$ weakly converge to zero, i. e.
 $$\tau(x_\alpha
a)\rightarrow0$$ for all $a\in
\mbox{Lin}(\bigcup\limits_{1<q\leq\infty}L^{q}(M, \tau))\cong
L^{\omega}(M, \tau)',$ where $\cong$ denote the isomorhism (see Proposition 2.10). By Proposition 3.4 $d$
is $t$-continuous, and hence  by [10, Proposition 8.6.5] $d$ is
weakly continuous. Thus $(d(x_\alpha)_{\alpha\in A}$ weakly
converges to zero, i. e.
$$\tau(d(x_\alpha)a)\rightarrow0\eqno (7)$$ for all $a\in
\mbox{Lin}(\bigcup\limits_{1<q\leq\infty}L^{q}(M, \tau)),$ and, in
particular, for all $a\in L^{1}(M,
\tau)\bigcap(\mbox{Lin}(\bigcup\limits_{1<q\leq2}L^{q}(M, \tau))).$

 Since $(L^{1}(M,
\tau)\bigcap(\mbox{Lin}(\bigcup\limits_{1<q\leq2}L^{q}(M,
\tau))))'\cong M+L^{\omega}_2(M, \tau),$   from  (7) it follows that
$(d(x_\alpha ))_{\alpha\in A}$ is *-weakly converging to zero in $M+
L^{\omega}_2(M, \tau).$ The proof is complete.  $\blacksquare$

As it was mentioned in Proposition  2.6 the Arens algebra
$L^{\omega}(M, \tau)$ is an ideal in $M+L^{\omega}_2(M, \tau),$ and
therefore any element $a\in M+L^{\omega}_2(M, \tau)$ generates a
spatial derivation on $L^{\omega}(M, \tau)$ defined as
$$d(x)=ax-xa, \quad x\in L^{\omega}(M, \tau).$$

In this connection a natural problem arises whether the converse
assertion is also true, i. e. can any derivation on the Arens algebra
be represented in this form?

The main result of the present work is the following theorem, which
answers this question in affirmative  and gives a complete
description of derivations on the Arens algebra $L^{\omega}(M,
\tau).$

{\bf Theorem 3.7.} \emph{Let   $M$ be a  von Neumann algebra with a
faithful normal semi-finite trace $\tau.$ Then any derivation $d$ on
  $L^{\omega}(M, \tau)$ is spatial, morever it is implemented by an element of $M+L^{\omega}_2(M, \tau),$ i. e.  }
$$d(x)=ax-xa, \quad x\in L^{\omega}(M, \tau)$$
  \emph{for some} $a\in M+L^{\omega}_2(M, \tau).$

Proof. Since  the trace $\tau$ is semi-finite  there exists a net of
projections $(e_\alpha)$ with $\tau(e_\alpha)<\infty$ for all
$\alpha,$ such that $e_\alpha\uparrow \textbf{1}.$ Consider the
derivations
$$d_{e_{\alpha}}: L^{\omega}(e_{\alpha}Me_\alpha, \tau_{e_{\alpha}})\rightarrow
L^{\omega}(e_\alpha Me_\alpha, \tau_{e_{\alpha}})$$ defined  as in
(4). By Corollary 3.3 there exist $a_\alpha\in
 L^{\omega}(e_{\alpha}Me_\alpha, \tau_{e_{\alpha}})$ such that
$$d_{e_{\alpha}} (x)=a_\alpha x-xa_\alpha, \quad x\in L^{\omega}(e_{\alpha}Me_\alpha, \tau_{e_{\alpha}}),$$
and moreover the net $(a_\alpha)$ is bounded in $M+ L^{\omega}_2(M,
\tau).$

By Corollary 2.11 the net $(a_\alpha)$ contains a subnet which
*-weakly converges in $M+ L^{\omega}_2(M, \tau).$ Without loss of
generality we may assume that $a_\alpha\rightarrow a$ for some $a\in
M+ L^{\omega}_2(M, \tau).$

Let  $x\in M\cap L^{1}(M, \tau).$ If $\alpha\geq\beta,$ then
$e_{\alpha}d(e_{\beta}xe_{\beta})e_{\alpha}=d_{e_{\alpha}}
(e_{\beta}xe_{\beta})=a_\alpha
e_{\beta}xe_{\beta}-e_{\beta}xe_{\beta}a_\alpha.$ Therefore
$$d(e_{\beta}xe_{\beta})=a
e_{\beta}xe_{\beta}-e_{\beta}xe_{\beta}a.\eqno (8)$$
 By Proposition 3.6 the derivation $d$  maps any weakly convergent net into a
 *-weakly
convergent  one. Therefore from (8) it follows that $d(x)=ax-xa$ for
all $ x\in M\cap L^{1}(M, \tau).$

Now since the set  $M\cap L^{1}(M, \tau)$ is  $t$-dense in
$L^{\omega}(M, \tau)$ from the $t$-continity of $d$ it follows that
$d(x)=ax-xa$ for all $ x\in L^{\omega}(M, \tau).$ The proof is complete. $\blacksquare$

  The following theorem can be proved in a way similar to the
  proof of the Theorem 3.7.

 {\bf Theorem 3.8.} \emph{Let   $M$ be a  von
Neumann algebra  with a faithful normal semi-finite trace $\tau.$
Then any derivation  $d$ on
  $L^{\omega}_2(M, \tau)$ is spatial, morever it is implemented by an element of $M+L^{\omega}_2(M, \tau),$ i. e.}
$$d(x)=ax-xa, \quad x\in L^{\omega}(M, \tau)$$
  \emph{for some} $a\in M+L^{\omega}_2(M, \tau).$

\textbf{Remark 4.} For any integer  $s\geq3$ put $L^{\omega}_s(M,
\tau)=\bigcap\limits_{p\geq s}^{\infty}L^{p}(M, \tau).$ It is not
difficult to  show that the space  $L^{\omega}_s(M, \tau)$ is also a
complete metrizable locally convex *-algebra, with the topology
generated by the family of norms $\{\|\cdot\|_p\}_{p\geq s},$ and
$L^{\omega}_s(M, \tau)$ is an ideal in $M+L^{\omega}_2(M, \tau).$
Moreover, one can also prove similarly to Theorem 3.7 that any
derivation of the algebra $L^{\omega}_s(M, \tau)$ is spatial and
implemented by an element of the algebra $M+L^{\omega}_2(M, \tau).$

\vspace{1cm}

\textbf{Acknowledgments.} \emph{The second and third named authors
would like to acknowledge the hospitality of the $\,$ "Institut
f\"{u}r Angewandte Mathematik",$\,$ Universit\"{a}t Bonn (Germany).
This work is supported in part by the DFG 436 USB 113/10/0-1 project
(Germany) and the Fundamental Research  Foundation of the Uzbekistan Academy of Sciences.}

\end{document}